\documentclass[11pt]{article}
\usepackage{latexsym,euscript,amsmath,amssymb,amsbsy,amsfonts,amsthm,amsopn,amstext,amsxtra,amscd}
\usepackage{epsfig}
\usepackage{graphics}
\usepackage{url}

\usepackage[english]{babel}

\theoremstyle{plain}
\newtheorem{theorem}{Theorem}

\theoremstyle{definition}

\newtheorem{point}{}
\newcommand{\N}{{\mathbb N}}

\newcommand{\D}{\text{$\mathcal{D}$}}

\begin{document}
\centerline{\bf On Buck's measurability of certain sets}
\centerline{Miilan Pa\v st\'eka}
\vskip0,5cm
{\bf Abstract.}{\it In the first part we construct some Buck measurable sets. In the second part we apply the Niven theorem for Buck's measure density to certain sets.}
\subsection{Notation}
$N$ the set of natural number, \newline
In accordance with algebra we shall use the following symbols: \newline
$r+(m)=\{n \in \N; n \equiv r \pmod{m}$\newline
$(m)=0+(m)$ \newline
$aS=\{as; s \in S\}$.
\subsection{Buck's measure density} This set function was firstly defined in 1946 by R. C. Buck in the paper \cite{BUC}.
The value
$$
\mu^\ast(S)=\inf\Big\{\sum_{i=1}^k \frac{1}{m_i}; S \subset \bigcup_{i=1}^k r_i+(m_i) \Big\}
$$
for $S \subset \N$ is called {\it Buck's measure density} of $S$. This set will be called {\it Buck measurable} if
$\mu^\ast(S)+\mu^\ast(\N \setminus S)=1$. The system of all Buck measurable sets we denote $\D_\mu$. This system is an algebra of sets. The restriction
$\mu = \mu^\ast|_{\D_\mu}$ is a finitely probability measure on this algebra. In the work \cite{BUC} is for each $\alpha \in [0,1]$ constructed a set $B_\alpha \in \D_\mu$ such that
$\mu(B_\alpha)=\alpha$, (see also \cite{h}, \cite{st}). This construction was later
In the paper \cite{PAS},(see also \cite{p1}), the following is proven:
\begin{theorem} If $A_1, A_2, A_3, \dots $ are such disjoint sets form $\D_\mu$ that
\begin{equation}
\label{weaksigma}
\lim_{N\to \infty} \mu^\ast\Big(\bigcup_{k=N}^\infty A_k\Big)=0,
\end{equation}
then the set $A=\cup_{k=1}^\infty A_k$ belongs to $\D$ and
$$
\mu(A)=\sum_{k=1}^\infty \mu(A_k).
$$
\end{theorem}
This leads to the following
\begin{point}
\label{1} Suppose that $b_i, i=1,2,3, $ is such increasing sequence of natural numbers that $b_i|b_{i+1}, i\in \N$.
Let $H_i; i=1,2,3, \dots $ be the such sets from $\D_\mu$ that
all elements from union of these sets are relatively prime with all $b_i, i\in \N$. Then the union
$H = \cup_{i=1}^\infty b_iH_i $ is Buck measurable and
$$
\mu(H) = \sum_{i=1}^\infty \frac{\mu(H_i)}{b_i}.
$$
\end{point}
This leads to shorter construction of the set $B_\alpha$ for $\alpha \in [0.1]$. We can suppose that $\alpha <1$. Thus this number has diadic expansion $\alpha=0,a_1a_2a_3 \dots$.
Let $n_1 < n_2 < \dots $ be the sequence of all such $n$ that
$a_n\neq 0$. Thus
$$
\alpha= \sum_{k}\frac{1}{2^{n_k}}.
$$
Let ${\bf O}$ be the set of all odd numbers. It holds
${\bf O} \in \D_\mu$ and $\mu({\bf O})=\frac{1}{2}$. Put
$B_\alpha =\cup_{k} 2^{n_k-1}{\bf O}$. Then $B_\alpha \in \D_\mu$ and
$$
\mu(B_\alpha)=\sum_k\frac{\mu({\bf O})}{2^{n_k-1}}=
\sum_k\frac{1}{2^{n_k}}=\alpha.
$$
\begin{point}Let $p$ be prime and $E=\{e_1 < e_2 < e_n < \dots \}$ be an increasing sequence of natural numbers. Denote
$N(p,E)$ the set of natural numbers containing $p$ in canonical representation only with the exponents from $E$.
It holds
$$
N(p,E)= \bigcup_{n=1}^\infty p^{e_n}(\N \setminus (p)).
$$
We see that {\bf \ref{1}} implies that $N(p, E)$ is Buck measurable and
$$
\mu(N(p,E))=\Big(1-\frac{1}{p}\Big)\sum_{n=1}^\infty \frac{1}{p^{e_n}}.
$$
\end{point}
\begin{point} This result can be generalized in the following way. Let primes $p_1 < p_2< \dots <p_k$ be given with the infinte sets of natural numbers $E_1, \dots, E_k$. Denote
$N(p_1, \dots, p_k, E_1, \dots , E_k):=N$ the set of all natural numbers containing $p_i$ in canonical representation with exponents from $E_i$, $i=1, \dots, k$. Then $N$ is Buck measurable and
$$
\mu(N)=
\prod_{i=1}^k\Big(1-\frac{1}{p_i}\Big)\prod_{i=1}^k\sum_{n_i\in  E_i } \frac{1}{p_i^{n_i}}.
$$
\end{point}
\subsection{Reminder systems}
Denote $R(S:m)=|\{s \pmod{m}; s \in S \}|$ where $S \subset \N$ and $m \in \N$. Suppose that $\{B_N\}$ is such sequence that for each $d \in \N$ such $N_0 \in \N$ exists that
$d|B_N$ for $N>N_0$.
In the paper \cite{PAS} the is proven that following equality
\begin{equation}
\label{lim}
\mu^\ast(S)= \lim_{N \to \infty}\frac{R(S:B_N)}{B_N}
\end{equation}
holds for each $S\subset \N$. \footnote{This says also that
$\mu^\ast(S)=P(cl(S))$, where the closure is considered
in the compact ring of polyadic integers and $P$ is a Haar measure on this ring, (see \cite{n}).} Ralph Alexander proved certain result concerning of union of sets for asymptotic density, (see \cite{a1}, \cite{p2}). Using (\ref{lim}) an analogy of this result can br proven for Buck's measure density:
\begin{point} Let $A_n, n=1,2,3, \dots $ be disjoint sets belonging to $\D_\mu$. Suppose that such convergent series with positive summands $\sum_{n=1}^\infty c_n$ exisits that for each $n, N \in \N$ the inequality
$$
\frac{R(A_n:B_N)}{B_N}\le c_n
$$
holds. The set $A=\cup_{n=1}^\infty A_n$ is Buck measurable and
$$
\mu(A)=\sum_{n=1}^\infty \mu(A_n).
$$
\end{point}
\subsection{Sets of zero Buck's measure density}
Ivan Niven proved in 1951 the result which characterise the sets of asymptotic density $0$ from "small" parts of given set, (see \cite{niv}). This result was later proved for Busk's measure density also (see \cite{p}).

Let $S \subset \N$ a $p$ prime. Denote
$S_p=\{s \in S; p|s \land p^2 \nmid s\}$.
\begin{theorem}
\label{niv}
Suppose that $\{p_i\}$ is such sequence of primes that
\begin{equation}
\label{divergence}
\sum_{i=1}^\infty \frac{1}{p_i} = \infty.
\end{equation}
Then for $S \subset \N$ we have $\mu(S)=0 \Leftrightarrow  \forall i=1,2,\dots \mu(S_{p_i})=0$.
\end{theorem}
\begin{point}Let the sequence of primes $\{p_n\}$ fulfils the condition (\ref{divergence}). If for the set $S$ the condition
$$\
\forall s \in S \forall n; p_n|s \Rightarrow p_n^2|s.
$$
Theorem \ref{niv} yields $S$ is Buck measurable and $\mu(S)=0$.
\end{point}
Denote by $R_t$ the set natural numbers containing at most $t$ primes from $\{p_n; n=1,2,3, \dots\}$ with odd exponent in canonical representation. We prove
\begin{point}
\label{tttt} The set $R_t, t\in \N$ is Buck's measurable
and $\mu(R_t)=0$.
\end{point}

{\bf Proof.} If $t=0$ then $(R_0)_{p_i}= \emptyset, i=1,2,3, \dots$. Thus $R_0 \in \D_{\mu}$ and $\mu(R_0)=0$.

Suppose now that $R_{t-1}\in \D_{\mu}$ and $\mu(R_{t-1})=0$.
For the set $R_t$ we have $(R_t)_{p_i} \subset p_iR_{t-1}$.
This yields $R_t \in \D_\mu$ and $\mu(R_t)=0$. \qed

This implies (see also \cite{lt1}, \cite{h}):
\begin{point} Let $P_t$ be the set of natural number containing at most $t$ primes in canonical representation.
Then $P_t \in \D_\mu$ and $\mu(P_t)=0$ for $t\in \N$.
\end{point}

Let $\tau(n)$ be a number of divisors of given $n \in \N$.
This function can be represented by canonical decomposition in the form
\begin{equation}
\label{exponent}
n=p_1^{\alpha_1}\dots _k^{\alpha_k} \Rightarrow
\tau(n)= (\alpha_1+1)\dots (\alpha_k+1).
\end{equation}
Viliam Fur\'\i k, (see \cite{Fu}), was interested in the set
$$
R = \{n \in \N; \tau(n)|n\}.
$$
\begin{point} The set $R$ is Buck measurable and $\mu(R)=0$.
\end{point}

{\bf Proof.} Let us denote by $P_s$ set of naturals numbers
containing at most $s$ prime numbers with odd exponents in canonical decomposition. From {\bf{\ref{tttt}}} we get $P_s$ is Buck measurable and $\mu(P_s)=0$. The set mentioned above we can decompose
\begin{equation}
\label{decomp}
R=(R\cap P_s) \cup (R\cap(\N \setminus P_s)).
\end{equation}
Since $\mu(P_s)=0$ we get $\mu((R\cap P_s))=0$.
If $n \in R\cap(\N \setminus P_s)$ then $\tau(n)|n$. The number $n$ contains at least $s+1$ primes with odd exponents in canonical representation, thus from (\ref{exponent}) we have $2^{s+1}| n$. This yields
$R\cap(\N \setminus P_s)\subset (2^{s+1})$ and so taking account (\ref{decomp}) we get
$$
\mu^\ast(R)\le \frac{1}{2^{s+1}}.
$$
Considering $s \to \infty$ we can conclude that $R$ is Buck measurable and $\mu(R)=0$. \qed

\end{document}